\documentclass[12pt]{article}
\usepackage{amssymb,amsthm,amsmath,amsfonts,latexsym,tikz,hyperref,ytableau,authblk,stackengine}
\stackMath
\usepackage[hmargin=1in,vmargin=1in]{geometry}

\newtheorem{thm}{Theorem}[section]
\newtheorem{prop}[thm]{Proposition}
\newtheorem{cor}[thm]{Corollary}
\newtheorem{lem}[thm]{Lemma}
\newtheorem{conj}[thm]{Conjecture}
\newtheorem{exa}[thm]{Example}

\DeclareMathOperator{\syt}{SYT}

\newcommand{\vb}{{\ol v}}
\newcommand{\phit}{\tilde{\phi}}
\newcommand{\Vh}{\hat{V}}
\newcommand{\vh}{\hat{v}}

\newcommand{\ben}{\begin{enumerate}}
\newcommand{\een}{\end{enumerate}}
\newcommand{\ble}{\begin{lem}}
\newcommand{\ele}{\end{lem}}
\newcommand{\bth}{\begin{thm}}
\renewcommand{\eth}{\end{thm}}
\newcommand{\bpr}{\begin{prop}}
\newcommand{\epr}{\end{prop}}
\newcommand{\bco}{\begin{cor}}
\newcommand{\eco}{\end{cor}}
\newcommand{\bcon}{\begin{conj}}
\newcommand{\econ}{\end{conj}}
\newcommand{\bde}{\begin{defn}}
\newcommand{\ede}{\end{defn}}
\newcommand{\bex}{\begin{exa}}
\newcommand{\eex}{\end{exa}}
\newcommand{\barr}{\begin{array}}
\newcommand{\earr}{\end{array}}
\newcommand{\btab}{\begin{tabular}}
\newcommand{\etab}{\end{tabular}}
\newcommand{\beq}{\begin{equation}}
\newcommand{\eeq}{\end{equation}}
\newcommand{\bea}{\begin{eqnarray*}}
\newcommand{\eea}{\end{eqnarray*}}
\newcommand{\bal}{\begin{align*}}
\newcommand{\bce}{\begin{center}}
\newcommand{\ece}{\end{center}}
\newcommand{\bpi}{\begin{picture}}
\newcommand{\epi}{\end{picture}}
\newcommand{\bpp}{\begin{picture}}
\newcommand{\epp}{\end{picture}}
\newcommand{\bfi}{\begin{figure} \begin{center}}
\newcommand{\efi}{\end{center} \end{figure}}
\newcommand{\bprf}{\begin{proof}}
\newcommand{\eprf}{\end{proof}\medskip}

\newcommand{\bsl}{\begin{slide}{}}
\newcommand{\esl}{\end{slide}}
\newcommand{\bfr}{\begin{frame}}
\newcommand{\efr}{\end{frame}}

\newcommand{\hqed}{\hfill \qed}

\newcommand{\eqqed}[1]{$\rule{1ex}{0ex}\hfill{\dil#1}\hfill\qed$}

\newcommand{\ol}{\overline}

\newcommand{\hs}[1]{\hspace{#1}}
\newcommand{\hso}[1]{\hspace{-1pt}}
\newcommand{\vs}[1]{\vspace{#1}}
\newcommand{\qmq}[1]{\quad\mbox{#1}\quad}

\newcommand{\emp}{\emptyset}

\newcommand{\sbe}{\subseteq}

\newcommand{\setm}{\setminus}

\newcommand{\Cong}{\equiv}

\newcommand{\Ph}{\hat{P}}
\newcommand{\ph}{\hat{p}}

\newcommand{\ptn}{\vdash}

\newcommand{\case}[4]{\left\{\barr{ll}#1&\mbox{#2}\\#3&\mbox{#4}\earr\right.}

\newcommand{\ce}[1]{\lceil #1 \rceil}

\def\<{\langle}
\def\>{\rangle}

\newcommand{\ra}{\rightarrow}

\newcommand{\io}{\iota}
\newcommand{\ka}{\kappa}
\newcommand{\la}{\lambda}
\newcommand{\mut}{\tilde{\mu}}
\newcommand{\om}{\omega}

\newcommand{\La}{\Lambda}

\newcommand{\bbN}{{\mathbb N}}
\newcommand{\bbP}{{\mathbb P}}

\newcommand{\bbR}{{\mathbb R}}

\newcommand{\bbZ}{{\mathbb Z}}
\newcommand{\cA}{{\cal A}}

\newcommand{\cB}{{\cal B}}

\newcommand{\cD}{{\cal D}}

\newcommand{\cE}{{\cal E}}

\newcommand{\cM}{{\cal M}}

\newcommand{\cO}{{\cal O}}
\newcommand{\cP}{{\cal P}}

\newcommand{\cR}{{\cal R}}

\newcommand{\cT}{{\cal T}}

\newcommand{\Fb}{\ol{F}}

\newcommand{\pb}{\ol{p}}
\newcommand{\Pb}{\ol{P}}

\DeclareMathOperator{\Mod}{mod}

\newcommand{\dil}{\displaystyle}

\begin{document}
\pagestyle{plain}

\title{Partitions with parity restrictions: a bijective approach
}
\author[1]{William Keith}
\affil[1]{Department of Mathematics, Michigan Technological University, Houghton, MI 49931}
\author[2]{Bruce E. Sagan}
\affil[2]{Department of Mathematics, Michigan State University, East Lansing, MI 48824}

\date{\today\\[10pt]
	\begin{flushleft}  
	                                       \\[5pt]
	\end{flushleft}}

\maketitle

\begin{abstract}
There has been recent interest in integer partitions whose parts satisfy parity restrictions: for example, those where  all the odd parts are distinct, or those where all the even parts are larger than the odd parts.  Often results about such partitions have been obtained by algebraic manipulation of generating functions.  We show that a number of these identities can be proved in a bijective, and sometimes simpler, manner.
\end{abstract}

\section{Introduction}

We will write $\bbZ$ and $\bbN$ for the integers and nonnegative integers, respectively.  If $n\in\bbN$ then let
$$
[n]=\{1,2,\ldots,n\}.
$$
For a finite set $S$ we will use the notations $\#S$ or $|S|$ to denote the cardinality of $S$.  For our generating functions in the variable $q$ we will use the {\em Pochhammer symbol}
$$
(a;q)_n = \prod_{i=0}^{n-1} (1-a q^i)
$$
for $n\in\bbN\uplus\{\infty\}$ where $\uplus$ is disjoint union.

We will now review some of the necessary definitions concerning partitions.  More information can be found in the book of Andrews~\cite{and:tp}.
A {\em partition} of $n$ is a weakly decreasing sequence of positive integers called {\em parts}.
$$
\la=(\la_1,\la_2,\ldots,\la_l)
$$
such that 
$|\la|:=\sum_i\la_i = n$.  
In this case we will also write $\la\ptn n$ or $|\la|=n$ where $|\la|=\sum_i\la_i$ is called the {\em size} of $\la$.  The {\em length} of $\la$ is the number of parts and we let
$$
\ell(\la) = \text{ length of $\la$}.
$$
When convenient we will extend $\la$ by letting $\la_i=0$ if $i>l$ and we will do this when necessary without further comment.

Let
$$
P(n) =\{\la \mid \la\ptn n\}
$$
and
$$
p(n) =\# P(n).
$$
In general, we will use upper case letters for sets and the corresponding lower case letters for their cardinalities.  For the set of all partitions, we use the notation
\beq
\label{cP}
\cP=\biguplus_{n\ge0} P(n).
\eeq
The same convention of dropping $n$ to not restrict the size of a family of partitions will be used in the sequel.

The {\em length} of $\la$, denoted $\ell(\la)$, is the number of parts $\la_i$.
  Using this notation, $\ell(\la)$ is the number of nonzero parts.  We also use {\em multiplicity notation}, writing $\la$ as a multiset
$$
\la=\{\{1^{m_1}, 2^{m_2},\ldots, n^{m_n}\}\}
$$
where $m_i$ is the number of times that $i$ appears in $\la$.  We will sometimes simplify this notation by leaving out the set braces and commas, as well as possibly writing the parts in a different order.  In a similar vein, given $\la$ we let
$$
m_i(\la) = \text{the multiplicity of $i$ as a part of $\la$}.
$$

The {\em Young} or {\em Ferrers diagram} of $\la$ is an array of left-juetified boxes (also called cells) with $\la_i$ boxes in row $i$ from the top.  We make no distinction between a partition and its Young diagram.  
We will index the boxes of $\la$ as in a matrix with $(i,j)$ being the cell in row $i$ and column $j$.
The {\em main diagonal} of $\la$ is 
$$
D(\la) = \{(1,1),\ (2,2),\ \ldots,\ (m,m)\}
$$
where $m$ is the largest index $i$ such that $(i,i)\in\la$ and is called the {\em length} of $D(\la)$.
The {\em conjugate} of $\la$, $\la^t$, is obtained by reflecting its Young diagram about the main diagonal.  Call $\la$ {\em self conjugate} if $\la^t=\la$.

Some of the results we consider will concern overpartitions.
An {\em overpartition} is a partition $\la$ where one is permitted to overline the first in any sequence of equal parts.  We ket
$$
\Pb(n) =\{\la \mid \text{$\la$ is an overpartition of $n$}\}.
$$
For example
$$
\Pb(3) =\{3,\ \ol{3},\ 2\, 1,\ \ol{2}\,1,\ 2\,\ol{1},\ \ol{2}\, \ol{1},\ 
1\, 1\, 1,\ \ol{1}\, 1\, 1\}.
$$
In general, putting a bar over the notation for a set  partitions represents considering all overpartitions whose underlying partition is in the original set.  
And overlining a cardinality is the cardinality of the corresponding overlined set.

There is a large literature about partitions with various parity restrictions on their parts.  Often, theorems about these partitions are obtained by manipulation of generating funcitons.  The purpose of the current work is to show how some of these results can be proved bijectively.  The rest of this paper is structured as follows.  In the next section we consider partitions of $n$ where all even parts are smaller than all odd parts, denoted $P_e^o(n)$, and those where the reverse is true, $P_o^e(n)$.  In particular, we give bijective proofs of various identities of Andrews~\cite{and:ipe,and:pps}.
Section~\ref{mr} is concerned with a subset of $P_e^o(n)$ obtained by making multiplicity restrictions on the various parts.  We use the lattice-path model for partitions to give combinatorial demonstrations of results of Chern~\cite{che:npe} and Passary~\cite{pas:spf}.  In Section~\ref{ro} we give a bijective proof of a theorem of Banerjee, Bringmann, and Dixit~\cite{BBD:roc} concerning overpartitions with certain parity restrictions.  Section~\ref{tom} is devoted to the study of the third order mock theta function.   Andrews~\cite{and:ipe} showed that its coefficients count certain parity-restricted partitions and we give a combinatorial proof of this fact along with bijections to various other partition models.   In addition to restricting the sizes of parts of different parity, one can also specify whether they are distinct or not.  Section~\ref{dvr} gives a bijective proof of an identity of Bringmann and Jennings-Shaffer~\cite{BJS:nap} about such partitions.  Guadalupe~\cite{gua:cal} proved a congruence modulo $2$ for the number of partition triples with certain parity conditions on their parts.  Our bijective demonstration in Section~\ref{pt}  generalizes this result to any prime modulus.
Section~\ref{prt} considers standard Young tableux having parity restrictions on their shapes.  We give bijections for proving results of Matsakis and Vendervelde~\cite{MV:mib} as well as Hemmer, Straub, and Westrem~\cite{HSW:eki} in this setting.  We end with a section of comments and open problems.

\section{Parity-separated parts}
\label{psp}

The purpose of this section is to study partitions with parts separated by parity.
A partition has {\em parts separated by parity} if all its even parts are smaller than all its odd parts or vice versa.  Both include the cases where all parts are even or all are odd.  We let $P_e^o(n)$, respectively $P_o^e(n)$, be the set of all partitions of $n$ with every even part is smaller than every odd part, respectively every odd part is smaller than every even part.  The cardinalities of these sets are denoted
$$
p_e^o(n) = \# P_e^o(n) \qmq{and} p_o^e(n) = \#P_o^e(n).
$$
For example
$$
P_o^e(5) = \{5,41,3 1^2,2^2 1, 2 1^3, 1^5\}
$$
so that $p_0^e=6$.  

To state the next result, we will need two other sets of partitions.  First of all, let
$$
P_{e,1}(n) = \{\la\ptn n \mid \text{every $\la_i$ is even or $1$}\}.
$$
Call a partition $\la$ {\em palindromic} if the number of multiplicities $m_i$ which are odd is at most one and let
$$
P_p(n) =  \{\la\ptn n \mid \text{$\la$ is palindromic}\}.
$$
These partitions are so called because in this case it is possible to arrange the parts of $\la$ in a palindrome.  For example
$$
P_p(5) = \{5, 3 1^2, 2^2 1, 1^5\}.
$$
The reader will note that the  bijection in the proof below  of 
$\#P_e^o(n) = \#P_{e,1}(n)$ is similar to the one  one defined by Franklin in his proof
of the Pentagonal Number Theorem~\cite[Theorem 6.1]{and:tp}.  A different Franklin-like bijection was applied to partitions with parts separated by parity by Fu and Tang~\cite{FT:pps}.

\bth[{\cite[Proposition 2.1]{and:ipe}}]
\label{p_e^o(n)}
For all $n\ge0$ we have
$$
\#P_e^o(n) = \#P_{e,1}(n) = \#P_p(n).
$$
\eth
\bprf
$\#P_e^o(n) = \#P_{e,1}(n)$.
Define a function $f:P_e^o(n) \ra \cP_{e,1}(n)$ as follows.
Suppose $\la=(\la_1,\ldots,\la_l)\in P_e^o(n)$ has odd parts 
$\la_1, \ldots,\la_k$.  We let $f(\la)=\mu$ where
$$
\mu=(\la_1-1,\ldots,\la_k-1,\la_{k+1},\ldots,\la_l,1^k).
$$
An example follows the proof.
Note that $f(\la)$ is well defined:  First of all,  $\mu$ is a partition since,  by the definition of $k$ and the parity constraints, we have $\la_k>\la_{k+1}$ so that $\la_k-1\ge \la_{k+1}$.  From the definition of $\mu$ we clearly have $|\mu|=|\la|=n$.  And  $\mu$ has only even parts except for $1$'s since one was subtracted from each odd part.  So, $\mu\in P_{e,1}(n)$.

To see that $f$ is a bijection, we construct its inverse.  Given 
$\mu\in P_{e,1}(n)$ having $k$ parts equal to one, we define $f^{-1}(\mu) = \la$ where
$$
\la = (\mu_1+1,\ldots,\mu_k+1,\mu_{k+1},\mu_{k+2}, \ldots)
$$
where, as usual, $\mu_i=0$ for $i>\ell(\la)$.
Checking that $f^{-1}$ is well defined and indeed the inverse of $f$ is straightforward so we omit the details.

\medskip

$\#P_{e,1}(n) = \#P_p(n)$.
We now define a function $g:\cP_{e,1}(n)\ra P_p(n)$.
Suppose $\mu\in \cP_{e,1}(n)$ has $k$ parts which are one and let $\mut$ be $\mu$ with those parts removed.  We let $g(\mu)=\nu$ where
$$
\nu=\mut^t\uplus(k).
$$
recalling that $t$ denotes conjugation and $\uplus$ is disjoint union of the two partitions viewed as multisets.  Again, an example will be found after the proof.
Since $\mut$ only has even parts, $\mut^t$ has all parts repeated twice.  It follows that $\nu$ has at most one part of odd multiplicity, namely $k$, and so is in $P_p(n)$.  The construction of the inverse is easy and left to the reader.
\eprf

\begin{figure}
    \centering
\ydiagram[*(gray)]
{6+1,6+1,4+1}
*[*(white)]{6,6,4,4,2}
\hs{20pt}
$\stackrel{f}{\mapsto}$
\hs{20pt}
\ydiagram[*(gray)]
{6+0,6+0,4+0,4+0,2+0,0+1,0+1,0+1}
*[*(white)]{6,6,4,4,2}

\vs{20pt}

\hspace*{130pt}
$\stackrel{g}{\mapsto}$
\hs{20pt}
\ydiagram[*(gray)]
{5+0,5+0,4+0,4+0,0+3,2+0,2+0}
*[*(white)]{5,5,4,4,0,2,2}
    \caption{The maps $f:P_e^o(n) \ra \cP_{e,1}(n)$ and $g:\cP_{e,1}(n)\ra P_p(n)$ of Theorem~\ref{p_e^o(n)}}
    \label{p_e^0(27)}
\end{figure}

Consider the partition $\la=(7,7,5,4,2)\in P_e^o(25)$ as shown on the left in Figure~\ref{p_e^0(27)}.  The gray boxes are the ones to be subtracted.  The image $\mu$ under $f$ is displayed on the right of the same line with the parts of size $1$ in grey.  Then applying $g$ one obtains the partition in the bottom line.

For our next result we will use the notation
\begin{align*}
P_e(n) &=\{\la\ptn n \mid \text{all parts of $\la$ are even}\},\\
P_o(n) &=\{\la\ptn n \mid \text{all parts of $\la$ are odd}\},\\
\Ph_o^e(n) &=\{\la\in P_o^e(n) \mid \text{$\la$ has no $1$'s and at least one odd part}\}.
\end{align*}
When the following statement is given in Andrews' paper~\cite{and:pps}, the condition of having at least one odd part is missing.   Andrews gives a proof of it which combines his anti-telescoping method~\cite{and:dpf} with a bijection.  Our proof will be purely bijective. 
Note that we will continue to use the notation $f$ for our bijection in the next result, but it is not the same function as defined previously.  This should cause no confusion because of the different context and we will also reuse letters for bijections in future results.
\bth[{\cite[Theorem 2.1]{and:pps}}] 
\label{hp_o^eThm}
We have
$$
\ph_o^e(n) = p_o(n) - p_e(n) -p_e(n-1).
$$
\eth
\bprf
It suffices to define a bijection
$$
f:P_o(n)\ra \Ph_o^e(n) \uplus P_e(n) \uplus P_e(n-1).
$$
We first consider the case when $n$ is even so that $P_e(n-1)=\emp$ so that we must define
$f:P_o(n)\ra \Ph_o^e(n) \uplus P_e(n)$.  An example corresponding to Figure~\ref{hp_o^eFig} follows the demonstration.

\begin{figure}
    \centering
 \ydiagram[*(gray)]
{3+0,3+0,3+0,3+0,1,1}
*[*(white)]{3,3,3,3}
\hs{20pt}
$\stackrel{f}{\mapsto}$ 
\hs{20pt}
 \ydiagram[*(gray)]
{3+1,3+1,3+0,3+0}
*[*(white)]{3,3,3,3}

\vs{20pt}

 \ydiagram[*(gray)]
{3+0,3+0,1+0,1,1,1}
*[*(white)]{3,3,1}
\hs{20pt}
$\stackrel{f}{\mapsto}$ 
\hs{20pt}
 \ydiagram[*(gray)]
{3+1,3+1,1+1}
*[*(white)]{4,4,2}
    \caption{The map $f:P_o(n)\ra \Ph_o^e(n) \uplus P_e(n) \uplus P_e(n-1)$  of Theorem~\ref{hp_o^eThm}}
    \label{hp_o^eFig}
\end{figure}

Since $n$ is even, any $\la=(\la_1,\la_2,\ldots) \in P_o(n)$ has even length. 
Let 
$$
k=\min(m_1(\la),\ell(\la)/2)
$$
and define $f(\la)=\mu$ where
\beq
\label{MuDef}
\mu=
\case{(\la_1+1,\ldots,\la_k+1,\la_{k+1},\ldots,\la_{l-k}, 0, 0,\ldots)}{if $m_1(\la)<\ell(\la)/2$,}
{(\la_1+1,\ldots,\la_k+1, 0, 0,\ldots)}{if $m_1(\la)\ge\ell(\la)/2$.}
\eeq
If $m_1(\la)<\ell(\la)/2$ then, by the way the $1$'s are added to the parts of $\la$, the partition $\mu$ will no longer have any $1$'s and its smallest part will be odd. Similar reasoning shows that all even parts will be greater than all odd parts.  
And clearly $|\mu|=|\la|=n$.
It follows that $\mu\in\Ph_o^e(n)$.  On the other hand, if
$m_1(\la)\ge\ell(\la)/2$ then every nonzero part of $\mu$ will be one more than a part of $\la$ and so even..  Thus $\mu\in P_e(n)$.

For $f^{-1}$, suppose $\mu=(\mu_1,\ldots,\mu_l)\in \Ph_o^e(n) \uplus P_e(n)$
and let $k$ be the number of even parts of $\mu$. We let $f^{-1}(\mu)=\la$ where
$$
\la=(\mu_1-1,\ldots,\mu_k-1,\mu_{k+1},\ldots,\mu_l,1^k)
$$
where, as usual, $1^k$ represents $k$ parts equal to $1$.  The fact that this is a well-defined and the inverse can be checked by the reader.

Finally, suppose that $n$ is odd so that $P_e(n)=\emp$ and we want to define a function
$f:P_o(n)\ra \Ph_o^e(n) \uplus P_e(n-1)$.
We keep the notation of the even case except with $\ell(\la)/2$ replaced by the ceiling $\ce{\ell(\la)/2}$ everywhere.  If $m_1(\la)<\ce{\ell(\la)/2}$ then $\mu=f(\la)$ is still defined  as in~\eqref{MuDef}.  On the other hand, if 
$m_1(\la)\ge\ce{\ell(\la)/2}$ then $\la$ has at least one $1$.
Let $\tilde{\la}$ be the result of removing a $1$ from $\la$ so that $\tilde{\la}\ptn n-1$ where $n-1$ is even.  So, we can let 
$f(\la):=f(\tilde{\la})$ where the right-hand side is defined by the even case.  The rest of the details are left as an exercise.
\eprf

Suppose  $\la=(3,3,3,3,1,1)$.  So,  $m_1(\la)= 2 < 3 =\ell(\la)/2$.  So $k=2$ and $f(\la)=(3+1,3+1,3,3)=(4,4,3,3)$.  See the top line in Figure~\ref{hp_o^eFig}.
On the other hand, suppose 
$\la=(3,3,1,1,1,1)$.  Now, $m_1(\la)=4>3=\ell(\la)/2$.  In this case $k=3$ and
$f(\la)=(3+1,3+1,1+1)=(4,4,2)$. This case is on the bottom in Figure~\ref{hp_o^eFig}.

\section{Multiplicity restrictions}
\label{mr}

We will now give a simple bijective proof of a result of Chern~\cite{che:npe}.
Let $V_e^o(n)$ be the set of all $\la\in P_e^o(n)$ satisfying the following two parity restrictions.
\ben
\item[(V1)]  All odd parts have even multiplicity.
\item[(V2)]  If even parts exist, then the largest has odd multiplicity and the rest even multiplicity.
\een
It will be convenient to let
$$
L(\la) = \text{ the largest even part of $\la$},
$$
where, by convention, we let $L(\la)=0$ if $\la$ has no even parts.  Now define
$$
V_e^o(n,m) = \{\la\in V_e^o(n) \mid L(\la)\Cong m\ (\Mod 4)\}.
$$
By (V1) and (V2) it is clear that the only possible values for $m$ above are $0$ and $2$.

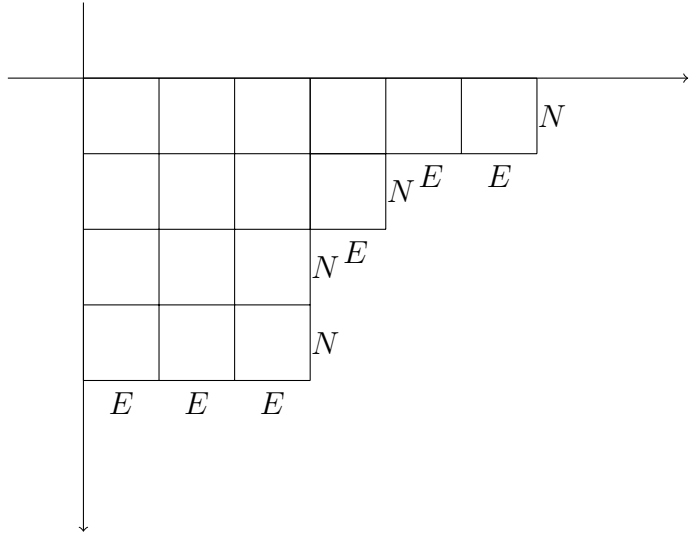
\begin{figure}
    \centering
 \begin{tikzpicture}
\draw[->] (-1,0)--(8,0);
\draw[->] (0,1)--(0,-6);
\draw (0,0) grid (3,-4);
\draw (3,0) grid (6,-1);
\draw (3,-1) grid (4,-2);
\draw (.5,-4.3) node{$E$};
\draw (1.5,-4.3) node{$E$};
\draw (2.5,-4.3) node{$E$};
\draw (3.2,-3.5) node{$N$};
\draw (3.2,-2.5) node{$N$};
\draw (3.6,-2.3) node{$E$};
\draw (4.2,-1.5) node{$N$};
\draw (4.6,-1.3) node{$E$};
\draw (5.5,-1.3) node{$E$};
\draw (6.2,-.5) node{$N$};
 \end{tikzpicture}
    \caption{The lattice path $P(6,4,3,3)$}
    \label{fig:P(la)}
\end{figure}

For our next bijection, it will  be useful to think of partitions in terms of lattice paths.  Consider the Young diagram of $\la=(\la_1,\ldots,\la_l)$.  See Figure~\ref{fig:P(la)} for an example of this when $\la=(6,4,3,3)$.  Embed the diagram of $\la$ in fourth quadrant of the plane $\bbR^2$ so that its northeast corner is at the origin and each square is a of unit size.  Now the southwest boundary of $\la$ defines a lattice path $P(\la)$ starting  at $(0,-l)$, ending at $(\la_1,0)$, and taking unit steps east, $E$, and north, $N$.  In our example,
$$
P(\la) = E E E N N E N E E N.
$$
For any word $w$ of $E$'s and $N$'s we denote its length (number of elements) by $\#w$.  The example has $\#P(\la) =10$.
A {\em run} in $P(\la)$ is a maximal factor (consecutive subword) which consists solely of $E$'s or solely of $N$'s.  These are called {\em $E$-runs} and {\em $N$-runs}, respectively. Our example has a sequence of $6$ runs which are: $EEE$, $NN$, $E$, $N$, $EE$, and $N$.  The {\em coordinate} of a run $R$ is the constant value of the $x$-coordinates of its steps if $R$ is an $N$-run, or of the $y$-coordinates of its steps if $R$ is an $E$-run.  Finally, for any word $w$ in $E$ and $N$ we let
$$
w^* = \text{ $w$ reversed and with $E$'s and $N$'s interchanged}.
$$
In our perennial example, 
$$
P(\la)^* = ENNENEENNN.
$$
We will often use the following easily-proved observations connecting $\la$ and $P(\la)$.
\ben
\item[(P1)]  If $R$ is an $N$-run of $P(\la)$ of length $r$ with coordinate $i$ then
$$
r = m_i(\la).
$$
\item[(P2)]  If $R$ is an $E$-run of $P(\la)$ of length $r$ with coordinate $-j$ then
$$
r= \la_j-\la_{j+1}
$$
where $\la_{j+1}=0$ if $j=l(\la)$.
\item[(P3)] We have
$$
P(\la^t) = P(\la)^*.
$$
\een

The $P(\la)$ for $\la\in V_e^o(n)$ have a particularly nice form.
\begin{lem}
\label{lem:P(la)}
We have $\la\in V_e^o(n)$ if and only if one of the two following conditions hold:
\ben
\item[(a)] $P(\la)$ has exactly one run of odd length which is the first or last, or
\item[(b)] $P(\la)$ has exactly two runs of odd length which are consecutive and begin with an $N$-run.
\een
\end{lem}
\bprf
We will only prove the forward implication as the reverse is similar.  We first show that $P(\la)$ has at most two runs of odd length.  For suppose it had three.  Then $P(\la)$ would have to contain two $N$-runs or two $E$-runs.  If the former holds then, by (P1), $\la$ has two parts of odd multiplicity which is impossible.  If the latter is true then, by (P2), the parts of $\la$ shift parity twice which is also a contradiction.

The proof now breaks into two parts giving items (a) and (b) in the statement of the lemma.  Part (a) corresponds to the $\la$ that have only even or only odd parts and the demonstration is similar to (b).  So we leave (a) to the reader.  As for (b), by the same reasoning as in the previous paragraph, one of the two runs must be an $N$-run and the other an $E$-run.  And, working from the smallest part of $\la$ to the largest, the even part of odd multiplicity occurs just before the change of parity of the parts.  This translates into the second half of the statement of (b).
\eprf

We are now in a position to give a bijective proof of an equality of Chern, answering his question to find a combinatorial proof.  We should note that Burson and Eichhorn~\cite{BE:pip} have also given a combinatorial proof of this result using a more complicated bijection.
\bth[{\cite[Theorem 1]{che:npe}}]
\label{n=2(mod4)}
If $n\Cong 2\ (\Mod 4)$ then
$$
v_e^o(n,0) = v_e^o(n,2).
$$
In particular, conjugation restricts to a bijection $V_e^o(n,0)\ra V_e^o(n,2)$.
\eth
\bprf
Obviously, it suffices to prove the second statement in the theorem.
From the previous lemma and (P3), we see that conjugation restricts to a bijection from $V_e^o(n)$ to itself for any $n$.  So, we only need to show that when $n\Cong 2\ (\Mod 4)$ we have $\la\in V_e^o(n,0)$ if and only if $\la^t\in V_e^o(n,2)$.  We will give details of the forward direction since the converse is similar.  And, as in the demonstration of the lemma, there are now two parts depending on whether conditions (a) or (b) hold, where we will only do the latter.

Suppose, towards a contradiction, that $L(\la)\Cong L(\la^t) \Cong 0\ (\Mod 4)$.  Break the Young diagram sitting in $\bbR^2$ up into rectangles as follows.  For each $E$-run before and including the unique $E$-run of odd length, draw a line with the same $y$-coordinate from the run to the $y$-axis.  Similarly, for each $S$-run after and including the unique $N$-run of odd length, draw a line with the same $y$-coordinate from the run to the $x$-axis.  It is now easy to check that the number of squares in each rectangle is divisible by $4$, either because both sides of the rectangle have even length, or because one of sides have length divisible by $4$.  Thus $n \Cong 0\ (\Mod 4)$ which is the desired contradiction.
\eprf

There is an analogue of the previous result for overpartitions.
Again, the following result was derived by Chern using manipulation of power series and he asked for a combinatorial proof.
\bth[{\cite[Theorem 2]{che:npe}}]
If $n\Cong 2\ (\Mod 4)$ then
$$
\vb_e^0(n,0) = \vb_e^o(n,2).
$$
\eth
\bprf
The number of overpartitions with underlying partition $\la$ is just the number of distinct parts of $\la$.  So, by the previous theorem, it suffices to show that the number of distinct parts of $\la$ is the same as the number of distinct parts of $\la^t$.  But the number of distinct parts of $\la$ is the number of $E$-runs since each such run gives a new part size.  And the number of $E$-runs equals the number of $N$-runs since the lattice path of $\la$ starts with an $E$-run and ends with an $N$-run.  Since conjugation just interchanges the $E$- and $N$-runs, we are done.
\eprf

\begin{figure}
    \centering
$\phi^2\ =\ \ydiagram{7,7,7,7,4,4,4}$
\hspace{50pt}
$\phit^1\ =\ \ydiagram{5,5,2,2,2}$
    \caption{The partitions $\phi^2$ and $\phit^1$}
    \label{F:fig}
\end{figure}

We will now give a totally bijective proof of a result of Passary~\cite{pas:spf} characterizing the parity of $v_e^o(n)$.  A partially combinatorial proof was given by Fu and Tang~\cite{FT:pps}.
Our proof will use an involution on the self-conjugate elements of $V_e^o(n)$.  

We first describe the fixed points of this map.  Given a partition $\la$ we define its {\em run-length sequence} to be $R(\la)=[r_1,r_2,\ldots,r_m]$ where $r_i$ is the length of the $i$th run of the lattice path $P(\la)$.  Using our usual example, we have 
$$
R(6,4,3,3) = [3,2,1,1,2,1].
$$
For $k\ge1$ we define the partition $\phi^k$ to be the one whose 
run-length sequence is
$$
R(\phi^k) = [2k,2k-1,2k-1,2k].
$$
The Young diagram of $\phi^2$ is shown on the left in Figure~\ref{F:fig}.  It is easy to see that $\phi^k$ is self-conjugate and 
\beq
\label{|Fk|}
|\phi^k|= 4k(3k-1).
\eeq
And it  follows from Lemma~\ref{lem:P(la)} (b) that $\phi^k\in V_e^o(4k(3k-1))$.  We will also need the partitions $\phit^k$ for $k\ge1$ defined by
$$
R(\phit^k) = [2k,2k+1,2k+1,2k].
$$
Note that 
\beq
\label{|Ftk}
|\phit^k|= 4k(3k+1)
\eeq
and that $\phit^k$ is a self-conjugate member of $V_e^o(4k(3k+1))$.
We have $\phit^1$ displayed in Figure~\ref{F:fig} on the right.

\begin{figure}
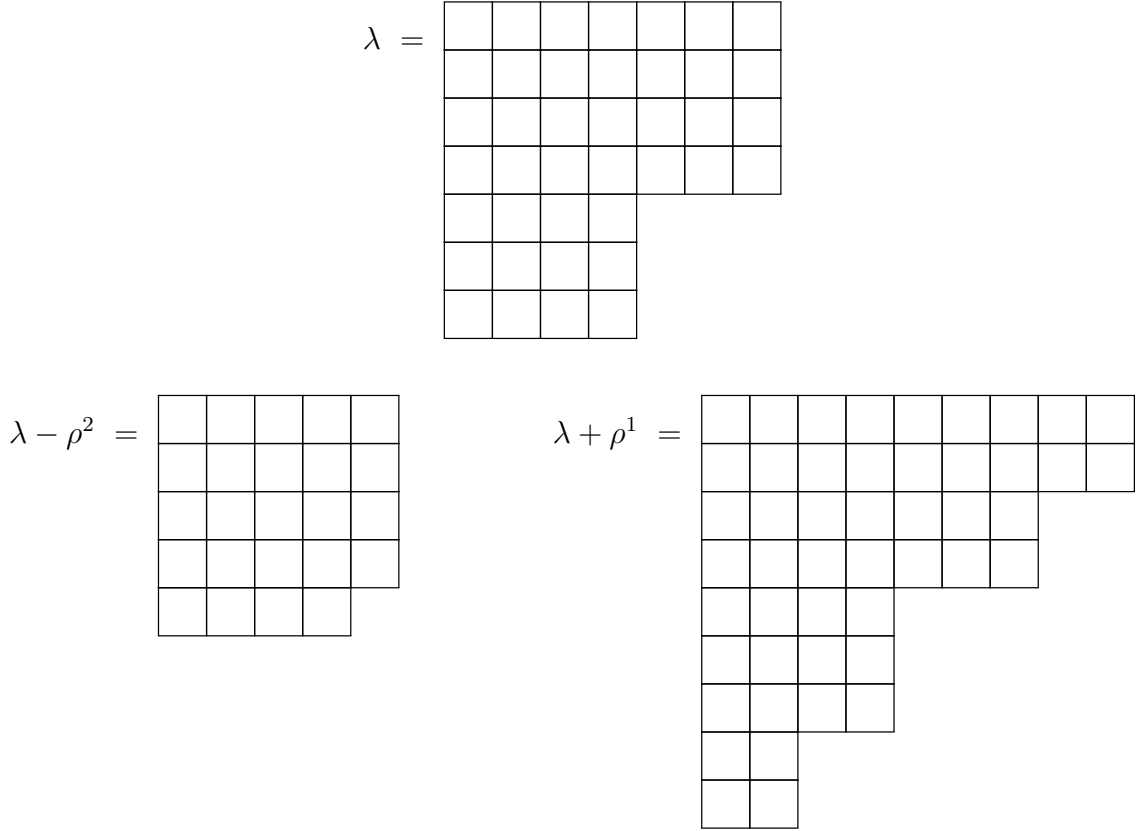

    \centering

\vspace*{-30pt}
    
  $\la\ =\ \ydiagram{7,7,7,7,4,4,4}$

  \vspace{20pt}

  $\la-\rho^2\ =\ \ydiagram{5,5,5,5,4}$
  \hspace{50pt}
  $\la+\rho^1\ =\ \ydiagram{9,9,7,7,4,4,4,2,2}$
    \caption{Subtracting and adding a double rectangle in $\la=(7,7,7,7,3,3,3)$}
    \label{+-rho:fig}
\end{figure}

To define the involution itself, we will need the following objects and operations.  For $k\ge1$, define the corrseponding {\em double rectangle} to be the pair of partitions 
$$
\rho^k = (\ \{\{2^{2k}\}\},\ \{\{(2k)^2\}\}\ ).
$$
Defining the size of a pair to be the sum of the sizes of its components, we have
$$
|\rho^k| = 8k.
$$
For $\la$  a partition whose last two columns and last two rows are all of length $2k$  we define
$$
\la-\rho^k = \text{ $\la$ with  its last two rows and last two columns removed}.
$$
See Figure~\ref{+-rho:fig} for the Young diagram of $\la-\rho^2$ where $\la=(7,7,7,7,3,3)$.  Similarly, if $\la$ is a partition satisfying the above conditions and $l\le k$, then we define
$$
\la+\rho^l = \text{ $\la$ with two final columns and two final rows added, all of length $2l$}.
$$
The preceding example is continued in  Figure~\ref{+-rho:fig} which also shows $\la+\rho^1$.  
Note that if $\la$ has last two rows and last two columns of length $k$ and $l\le k$ then
$$
\la-\rho^k+\rho^k = \la \text{ and } \la+\rho^l-\rho^l = \la
$$
where addition and subtraction are done left to right.

\begin{figure}
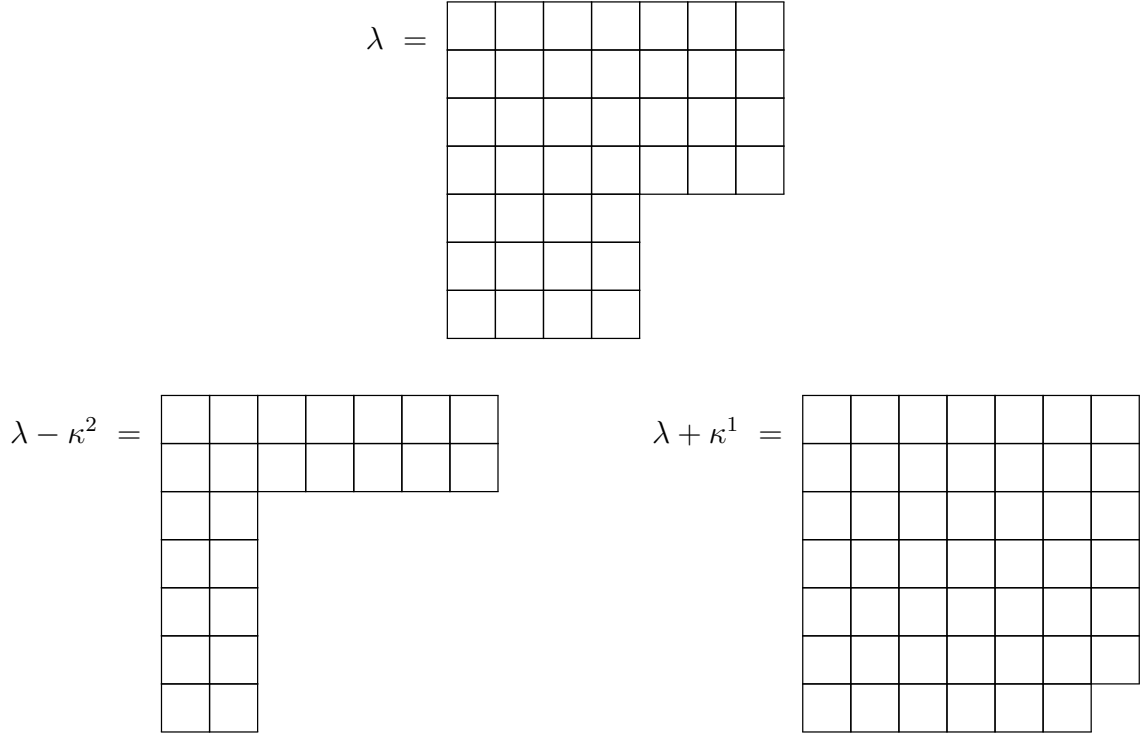

    \centering

\vspace*{-30pt}
    
  $\la\ =\ \ydiagram{7,7,7,7,4,4,4}$

  \vspace{20pt}

  $\la-\ka^2\ =\ \ydiagram{7,7,2,2,2,2,2}$
  \hspace{50pt}
  $\la+\ka^1\ =\ \ydiagram{7,7,7,7,7,7,6}$
    \caption{Subtracting and adding a double hook in $\la=(7,7,7,7,3,3,3)$}
    \label{+-ka:fig}
\end{figure}

We will also use  hooks.  The {\em hook} of $(i.j)\in\la$ is 
$$
H_{i,j} = \{(i,j')\in\la \mid j'\ge j\} \cup \{(i',j)\in\la \mid i'\ge i\}.
$$
For $k\ge1$, the {\em $k$th double hook} is the partition
$$
\ka^k = \{\{(2k+1)^2, 2^{2k-1}\}\}.
$$
Observe that
\beq
\label{ka:rho}
|\ka^k| = 8k = |\rho^k|.
\eeq
Suppose that $\la$ is a partition with main diagonal $M(\la)=\{(1,1),\ldots,(m,m)\}$ such that we have 
$H_{m-1,m-1}\uplus H_{m.m}=\ka^k$ for some $k$.  In that case we define
$$
\la-\ka^k = \text{ $\la$ with $H_{m-1,m-1}\uplus H_{m.m}$ removed}.
$$
An example is given in Figure~\ref{+-ka:fig}.  Now if $\la$ satifies the double hook subtraction conditions and $l<k$ then we define
$$
\la+\ka^l = \text{ $\la$ with boxes added so that $H_{m+1,m+1}\uplus H_{m+2,m+2}$ forms $\ka^l$}.
$$
Again, Figure~\ref{+-ka:fig} contains an example.  As with subtraction and addition of double rectangles, for $l<k$,
$$
\la-\ka^k+\ka^k = \la \text{ and } \la+\ka^l-\ka^l = \la.
$$

We now have everything in place to combinatorially prove the desired parity result.
\bth[\cite{pas:spf}]
For $n\ge0$ and $k\ge1$ we have
$$
v_e^o(n) \Cong \case{1\ (\Mod 2)}{if $n=4k(3k\pm1)$,}{0\ (\Mod 2)}{else.}
$$
\eth
\bprf
Consider the self-conjugate partitions in $V_e^o(n)$:
$$
\Vh_e^o(n) =\{\la\in V_e^o(n) \mid \la^t=\la\}.
$$
Since conjugation is an involution, all the cycles in its cycle decomposition as a permutation have length $1$ and $2$.  And, as we saw in the proof of Theorem~\ref{n=2(mod4)}, conjugation acts on $V_e^o(n)$.  So $v_e^o(n)$ and $\vh_e^o(n)$ have the same parity.  Thus, it suffices to find an involution $I$ on $\Vh_e^o(n)$ with exactly one fixed point for each $n=4k(3k\pm1)$.

Suppose that  $\la\in \Vh_e^o(n)$.  Then $\la$ satisfies Lemma~\ref{lem:P(la)} (b) since the partitions in (a) are not self-conjugate.  Let $2k$ be the length of the first run (and also the last run) of $P(\la)$ and $2l-1$  be the length of the two odd runs of $P(\la)$.  Define the map
$$
I(\la)=
\case{\la-\rho^k+\ka^k}{if $k<l$,}{\la-\ka^l+\rho^l}{if $k\ge l$,}
$$
where $I(\la)=\la$ if any of the addition or subtraction operations above are not defined.
We will show that $I$ is an involution whose fixed points are the $\phit^k$ when applied to the partitions where $k<l$.  The case $k\ge l$ is similar with fixed points the $\phi^k$.

We must first show that the addition and subtraction operations in the definition of $I$ yield valid partitions.   We first consider $\la-\rho^k$.  By definition of $k$, the first column and last row of $\la$ each contain $2k$ cells.  So, if the same is true of the penultimate column and row then subtraction of $\rho^k$ is possible.  This is clearly true if $P(\la)$ has at least six runs since then the second and penultimate run are of even length which is at least $2$.  So consider the case 
$R(\la)=[2k,2l-1,2l-1,2k]$ and suppose, towards a contradiction, that $2l-1=1$.  This implies that $l=1$.  But we are assuming that $k<l$ which cannot be.  Next we must check that $\ka^k$ can be added to $\la-\rho^k$.  Since $\ka^k$ has first row and first column of length $2k+1$, it suffices to show that the odd runs in $P(\la-\rho^k)$ are at least this long.  If $P(\la)$ has at least six runs, then the odd runs are not affected by subtracting $\rho^k$.  And $k<l$ implies $2k+1\le 2l-1$ so we are done in this case.  When there are four runs, we
have $R(\la-\rho^k)=[2k,2l-3,2l-3,2k]$ and the previous argument will still go through unless
 $k=l-1$.  But then $\la=\phit^k$ which is one of our desired fixed points.

 Next, we must show that $I$ is well defined in that $I(\la)\in\Vh_e^o(n)$.  It is easy to check that the operations of addition and subtraction preserve the properties in Lemma~\ref{lem:P(la)} and of being self-conjugate.  Also, using equation~\eqref{ka:rho},
 $$
 |\la-\rho^k+\ka^k| = |\la|-|\rho^k|+|\ka^k| = |\la| = n.
 $$
 Thus, $I(\la)$ has all the necessary properties to be in $\Vh_e^o(n)$.

 Lastly, we must check that $I$ is an involution.  Let
 $$
 \la'=\la-\rho^k+\ka^k
 $$
 and let $k'$ and $l'$ be the analogues for $\la'$ of $k$ and $l$ for $\la$.  Removing the last two rows and last two columns can only weakly increase $k$ so that $k'\ge k$.  And $l'$ is the length of the odd runs in $\ka^k$ so that $l'=k$.  It follows that $k'\ge l'$ so that $I(\la')$ is determined by the second case in the definition of $I$.  Using this and the fact that $l'=k$ gives
 $$
 I^2(\la) = I(\la-\rho^k+\ka^k) =(\la-\rho^k+\ka^k)-\ka^k+\rho^k=\la
 $$
 and the proof is complete.
 \eprf

\section{Restricted Overpartitions}
\label{ro}

We now give a bijective proof of an identity of 
Banerjee, Bringmann, and Dixit~\cite{BBD:roc}, answering a question posed by the authors.
Their result relates certain overpartitions with certain two-colored partitions, both having restrictions based on the parities of their parts.


Let $\overline{P}_{od}(n)$ be the set of overpartitions of $n$ with the following restrictions: 
\ben
\item[] \hspace*{-30pt} (OD1) Odd parts are distinct and overlined.
\item[] \hspace*{-30pt} (OD2) If the smallest part is even and occurs only once, it is not overlined.
\een
We note that in~\cite{BBD:roc} they consider overpartitions satisfying (OD2) and having distinct odd parts all of which are not overlined.  But clearly both sets of overpartitions have the same cardinality, and for our purposes the ones we have defined are more natural.

A {\em two-colored partition} is a partition $\la$ where each part $i$ has been colored red or blue denoted $i_r$ or $i_b$, respectively.  Note that the order of the colors does not matter so that, for example,
$$
(3_r,3_b,2_b,1_b,1_r,1_r) = (3_b,3_r,2_b,1_r,1_b,1_r).
$$
Let $G(n)$ be the set of two-colored partitions of $n$ such that
\ben
\item[(G1)] 
The smallest part size is odd, say $2k+1$ for some $k\ge0$, and
\item[(G2)] 
The only parts which can be colored blue are those of size $2k+2, 2k+4,\ldots,4k$.
\een

Finally, if $\la$ is a  usual partition with positive parts then we will denote by $\la0$ the partition obtained by appending a unique part equal to $0$ to $\la$.  To illustrate,
$$
(3,3,2,1,1)0 = (3,3,2,1,1,0).
$$
Define
$$
P_{ze}(n) = \{\la0 \mid \la\in P(n)\}
$$
Clearly $p_{ze}(n) = p(n)$.

\bth[{\cite[Equation (1.10)]{BBD:roc}}]
    For $n\ge1$ we have
    $$
 g(n) + \overline{p}_{od}(n-1) = 2 p(n-1).   
    $$
\end{thm}
\bprf
It suffices to find a bijection
$$
f:G(n)\uplus \Pb_{od}(n-1) \ra P_{ze}(n-1) \uplus P(n-1).
$$
An example of this map follows the proof.
If $\ka\in G(n)$ then define $f(\ka) = \la$ where $\la$ is the (uncolored) partition obtained from $\ka$ by the following two rules.
\ben
\item[(fr)] Replace a copy of the smallest part $(2k+1)_r$ in $\ka$ by $2k$ in $\la$, and replace all other red parts by (uncolored) copies of themselves.
\item[(fb)]  Replace each blue part $(2l)_b$ in $\ka$ by $l^2$ in $\la$.
\een
If $\om\in\Pb_{od}(n-1)$, then let $f(\om)=\mu$ where $\mu$ is the (not overlined) partition obtained from $\om$ as follows
\ben
\item[(fo)]  Replace each overlined part of $\om$ by a (not overlined) part of the same size in $\mu$.
\item[(fe)]  Replace each not overlined $2l$ in $\om$ by $l^2$ in $\mu$.
\een

To show that $f$ is bijective, first consider the subset $G'(n)\subseteq G(n)$ consisting of all $\ka$ with smallest part $1$.  So, $k=0$ in (G1) and (G2).  But this forces all parts of $\ka$ to be red.  It follows from (fr) that $f$ restricts to a bijection $f:G'(n)\ra P_{ze}(n-1)$.

Next consider $\ka\in G(n)\setm G'(n)$.  By (fr) and (fb), the blue parts of $\ka$ are the only ones which contribute part sizes in $\la$ in the range $[k+1,2k]$, except for the $(2k+1)_r$ which contributes a part $2k$.  This forces $2k$ to be the smallest part size of odd multiplicity.
Thus these $\la$ are all partitions of $n-1$ satisfying the following two properties.
\ben
\item[(i)] There is at least one even part of odd multiplicity.
\item[(ii)] if $2k$ is the smallest such part, then every part is larger than $k$.
\een

If $P'(n-1)$ is the set of partition in  $P(n-1)$ satisfying (i) and (ii), then we claim that the restriction $f:G'(n)\ra P'(n-1)$ is bijective.  To see this, we construct its inverse.  Given $\la\in P'(n-1)$ we define $\ka=f^{-1}(\la)$ as follows.
\ben
 \item[(f-r)] Replace one of the parts of size $2k$ with $(2k+1)_r$ and color all the parts of size at least $2k+1$ red.
 \item[(f-b)] Pair up the remaining parts in $[k+1,2k]$ and replace each pair $l^2$ with $(2l)_b$.
\een 
Since (f-r) and (f-b) are step-by-step inverses of (fr) and (fb), it is clear that the restriction under consideration is a bijection. 

Finally, let $P''(n)=P(n)\setm P'(n)$.  We will complete the proof if we can show that the restriction $f:\Pb_{od}(n-1)\ra P''(n-1)$ is  well-defined in that $f(\Pb_{od}(n-1))\sbe P''(n-1)$ and, in fact, is a bijection.  To show that $f$ is well-defined suppose, to the contrary, that $f(\om)=\mu$ where 
$\mu\in P'(n-1)$.  Note that by the definitions of $\Pb_{od}(n-1)$ and $f$, a part $l$ of $\mu$ has odd multiplicity if and only $\ol{l}$ is a part of $\om$.  So, by property (ii) of $P'(n-1)$, $\om$ contains a part $\ol{2k}$.  Furthermore, since all parts of $\mu$ smaller than $2k$ have even multiplicity and are bigger than $k$, the definition of $f$ shows that such parts correspond to halving parts of $\om$ which are larger than $2k$.  Thus $\om$ has a unique smallest part and it is of the form $\ol{2k}$, contradicting (OD2).

To see that $f:\Pb_{od}(n-1)\ra P''(n-1)$ is a bijection, we construct the inverse of this part of the restriction.  Given $\mu\in P''(n-1)$ we let $\om=f^{-1}(\mu)$ be defined by the following two conditions.
\ben
\item[(f-o)] For each part $l$ of $\mu$ with odd multiplicity, let $\ol{l}$ be a part of $\om$.
\item[(f-e)]  Pair up the remaining parts of $\mu$ and replace each pair $l^2$ with a non-overlined $2l$. 
\een
Again, this is a steep-by-step reversal of (fo)--(fe) and so this map is the inverse of the restriction as claimed.  This also completes the proof of the theorem.
\eprf

As an example of the map used in the previous theorem, suppose $\ka=(4_r,4_b,4_b,3_r,3_r,3_r)$.  Then a $3_r$ becomes a $2$ in $\la=f(\ka)$.  For the remaining red parts, we merely remove the color in $\la$.  
And finally the two $4_b$ parts are both halved to obtain
$$
\la= (4,3^2,2^5).
$$
On the other hand, if $\om=(\ol{4},\ol{2},2,2,\ol{1})$ then in $\mu=f(\om)$ the three overlined parts are kept without the overline.  And the two non-overlined $2$'s are both halved so that
$$
\mu=(4,2,1^5).
$$

\section{The third order mock theta function}
\label{tom}

The {\em third order mock theta function} is 
$$
\nu(q) = \sum_{k\ge0} \frac{q^{k^2+k}}{(-q;q^2)_{k+1}}.
$$
Define a sequence $\{p_\nu(n)\}_{n\ge0}$  by
$$
\nu(-q) = \sum_{n\ge 0} p_\nu(n) q^n.
$$
There are various known combinatorial interpretations of the $p_\nu(n)$.  One can be obtained directly from the generating function
\bpr[{\cite[equation (5.6)]{and:ipe}}]
\label{pnuThm}
For $n\ge0$, the coefficient $p_\nu(n)$ is the number of 
$\la=(\la_1,\ldots,\la_l)\ptn n$ satisfying
\ben
\item[(a)]  the even parts of $\la$ are distinct, and
\item[(b)] if $2m<\la_1$ for some $m\in\bbP$ then $2m$ is a part of $\la$.
\een
\epr
\bprf
The $k$th  summand in $\nu(-q)$ is $q^{k^2+k}/(q;q^2)_{k+1}$ since the power $k^2+k$ is even.  Write $k^2+k = 2+ 4 +\cdots+2k$.  So, this term counts partitions with even parts $2,4,\ldots,2k$ each appearing exactly once, and with odd parts $1,3,\ldots,2k+1$ each appearing any number of times (including zero).  This is equivalent to the description in the proposition.
\eprf

For our next bijections we need the notion of an odd Ferrers graph as introduced in~\cite{and:pds}.    A {\em Young tableau of shape $\la$}, $T$, is a filling of the Ferrers diagram $\la$ with positive integers such that the rows and columns weakly increase left to right and top to bottom. Let 
$$
T_{i,j} = \text{entry of $T$ in cell $(i,j)\in\la$}
$$
and
$$
|T|=\sum_{(i,j)\in\la} T_{i,j}.
$$
We extend the tableaux values outside of the shape $\la$ by letting
$$
T_{i,j} = 0 \text{ if $(i,j)\not\in\la$}.
$$
The {\em odd Ferrers graph} corresponding to $\la$ is the Young tableau $T^\la$ of shape $\la$ such that
$$
T^\la_{i,j}=\case{1}{if $i=1$ or $j=1$,}{2}{otherwise.}
$$
The tableau in Figure~\ref{fBA} is $T^\la$ for $\la=(7,4,4,3,1,1)$.
Define the sets
\begin{align*}
A(n) &=  \{\la\ptn n  \mid \text{$\la$ satisfies conditions (a) and (b) in Proposition~\ref{pnuThm}}\},\\
B(n) & = \{T^\la \mid \text{$\la$ is self conjugate and $|T^\la|=2n+1$}\},\\
C(n) &= \{\la\ptn 4n+1 \mid \text{$\la$ is self conjugate and all parts are odd}\}.
\end{align*}
The following result is mentioned in~\cite{and:ipe} based on statements in the Online Encyclopedia of Integer Sequences.
We note that a map similar to the one we use for proving $\#B(n)=\#C(n)$ was employed by Chern~\cite{che:pga} in giving a combinatorial demonstration of a result of Andrews~\cite{and:ipe}.

\begin{figure}
 \centering
 $
\begin{array}{rl}
     &  \\[85pt]
1+2/2   & = 2\\[3pt]
1+2+2/2 & = 4\\[3pt]
1+2+2+0/2 & = 5\\[3pt]
1+0+0+0+0/2&=1\\[3pt]
1+0+0+0+0+0/2&=1
\end{array}
$
\hs{10pt}
 \begin{ytableau}
     1&1&1&1&1&1&1\\
     1&2&2&2\\
     1&2&2&2\\
     1&2&2\\
     1\\
     1
 \end{ytableau}
 \hs{20pt}
$\stackrel{f}{\mapsto}$ 
\hs{20pt}
(5,4,2,1,1)
    
    \caption{The map $f:B(n)\ra A(n)$ of Theorem~\ref{ABCThm}}
    \label{fBA}
\end{figure}

\bth
\label{ABCThm}
For all $n\ge0$ we have
$$
\# A(n) = \# B(n) = \# C(n).
$$
\eth
\bprf
$\#A(n)=\#B(n)$.  We construct a bijection $f:B(n)\ra A(n)$ as follows.  We let $f(T^\la)=\mu$ where the parts of $\mu$ are the sums
\beq
\label{TlaSum}
T^\la_{i,1} + T^\la_{i,2}+\cdots+T^\la_{i,i-1}+(T^\la_{i,i}/2)
\eeq
for $2\le i\le \ell(\la)$.  An example of this map and the one for proving $\# B(n)=\# C(n)$ will be found following the proof with corresponding Figures~\ref{fBA} and~\ref{gBC}.

We need to show that $f$ is well defined in that $\mu$ satisfies conditions (a) and (b) in Theorem~\ref{pnuThm}.  Let $k$ be the length of the main diagonal of $\la$.  Then the sums~\eqref{TlaSum} for rows $2,3,\ldots,k$ are the even numbers $2,4,\ldots,2k-2$.  And the sums for rows $i>k$ are all odd and weakly decreasing as $i$ increases.  Furthermore, the maximum value of the sum for row $k+1$ is $2k-1$.  It is now easy to see that these conditions imply (a) and (b).

For the inverse map, we will suppose $\mu\in A(n)$ and construct $T^\la=f^{-1}(\mu)$.  It suffices to define the restriction of the shape $\la$ to the cells $(i,j)$ with $i\ge j$ since then there is a unique corresponding self-conjugate $\la$.  And once $\la$ is determined, then the entries of $T^\la$ are fixed.
Suppose that the even parts of $\mu$ are $2,4,\ldots,2k-2$. Then we let $(i,j)\in\la$ for all $j\le i\le k$.
Furthermore, for each odd part $2m-1$ of $\mu$  we add cells to $\la$ corresponding to a row of length $m$.  As usual, we leave the facts that $f^{-1}$ is well defined and the actual inverse of $f$ to the reader.

\medskip

$\#B(n)=\#C(n)$.  Note that there is a unique tableau $T^\la$ for any $\la$.  So we construct $g:B(n)\ra C(n)$ by ``blowing up" the squares of a self-conjugate partition $\la$ to form $\mu=f(T^\la)$ in the following way.  We will replace each cell of $T^\la$ other than $(1,1)$, which is invariant, by a set of cells which we will call a {\em tile}.  Specifically,  each $(i,j)\in\la/(1,1)$ is replaced by
\ben
\item[(T1)] a horizontal tile of dimensions  $1\times 2$ if $i=1$, or
\item[(T2)] a vertical tile of dimensions  $2\times 1$ if $j=1$, or
\item[(T3)] a square tile of dimensions  $2\times 2$ if $i,j\neq 1$.
\een

We must verify the $\mu\in C(n)$.  Note that for all $(i,j)\neq(1,1)$ we are replacing a cell labeled $T_{i,j}^\la$ with a tile containing $2T_{i,j}^\la$ cells.  So, $|T^\la|=2n+1$ implies that $|\mu|=4n+1$.  The fact that $\la$ is self-conjugate and the dimensions of the cells ensure that $\mu$ is self conjugate.  As far as all the parts of $\mu$ being odd, it is easy to check  that if $\la_i=k$ then $\mu_{2i-1}=\mu_{2i-2}=2k-1$ where $\mu_0$ should be ignored in the previous equality.
The fact that, for $\mu\in C(n)$, all parts except the first must occur with even multiplicity follows easily from properties (P1)--(P3) of the lattice path $P(\mu)$.

As far as $g^{-1}$, consider $\mu\in C(n)$. A simple induction on $|\mu|$  shows that $\mu$'s Young diagram can be decomposed into tiles satisfying (T1)--(T3) above.  To obtain $T^\la=g^{-1}(\mu)$, one replaces each tile in the first row or column by a cell labeled $1$, and all other tiles by  cells labeled $2$.
\eprf

\begin{figure}
    \centering
\begin{ytableau}
 1& *(lightgray) 1 & 1 & *(lightgray) 1 \\
*(lightgray) 1 & 2 & *(lightgray) 2 \\
 1 & *(lightgray) 2 \\
*(lightgray) 1 
\end{ytableau}
\hs{20pt}
$\stackrel{g}{\mapsto}$ 
\hs{20pt}
\begin{ytableau}
*(white)  & *(lightgray)  & *(lightgray)  &  &  & *(lightgray)  & *(lightgray) \\
*(lightgray)  &  &  & *(lightgray) & *(lightgray)\\
*(lightgray)  &  &  & *(lightgray) & *(lightgray)\\
 *(white) &  *(lightgray)  & *(lightgray) \\
 *(white) &  *(lightgray)  & *(lightgray) \\
*(lightgray) \\
*(lightgray) 
\end{ytableau}
    \caption{The map $g:B(n)\ra C(n)$ of Theorem~\ref{ABCThm}}
    \label{gBC}
\end{figure}

For the map $f:B(n)\ra A(n)$, suppose $T^\la$ is as given in Figure~\ref{fBA}.  Then the corresponding row sums are displayed to the left of $T^\la$.  And the resulting $\mu$ is gotten by sorting these values into the partition on the right.
An example of $g:B(n)\ra C(n)$ for $\la=(4,3,2,1)$ will be found in Figure~\ref{gBC}.  We have shaded the cells of $T^\la$ in a checkerboard fashion and done the same for the corresponding tiles to help clarify the construction.

We  note a that Chern~\cite{che:cpi} has defined a bijection between two different sets of objects which, when suitably modified, is the same as our function $f:B(n)\ra A(n)$.  Because the nature of Chern's map bears on a still-open search for a bijection regarding a fourth class associated with the previous theorem, we defer remarks on this to Subsection~\ref{mto}, where we discuss that problem.

\section{Distinct versus repeated parts}
\label{dvr}

We will consider the subsets of the partition sets we have considered thus far  by restricting either the even or the odd parts (or both) to be distinct.  This will be done in the notation by adding a $d$ or a $u$ depending on whether the parts with that parity are distinct or unrestricted.  For example
$$
P_{ed}^{ou} = \{\la \in P_e^o \mid \text{the even parts are distinct and the odd parts are unrestricted}\}.
$$
We apply the same rules to $P_o^e$.  

Andrews~\cite{and:pps} first derived expressions for the generating functions $F_{xy}^{zw}$ of the $p_{xy}^{zw}$.  New identities for three of them were given by  Bringmann and Jennings-Shaffer~\cite{BJS:nap}.  Here we give a combinatorial proof of one of them, slightly corrected to include partitions with all parts odd and distinct.
As usual, an example follows the theorem's proof.
\bth[{\cite[equation (1.1)]{BJS:nap}}]
\label{Fedod:Thm}
We have
$$
F_{ed}^{od}(q) = 
\frac{1}{1-q}(-q;q^2)_\infty - \frac{q}{1-q}(-q^2;q^2)_\infty.
$$
\eth
\bprf
Moving the negative term to the right-hand side we wish to show
\beq
\label{pedod}
F_{ed}^{od}(q) + \frac{q}{1-q}(-q^2;q^2)_\infty
=\frac{1}{1-q}(-q;q^2)_\infty.
\eeq
Consider pairs $(\la,s)$ where $\la$ is a partition, $s\in\bbN$, and 
$$
\case{\la\in P_{ed}^{od}}{if $s=0$,}
{\la\in P_{ed}}{if $s>0$.\rule{0pt}{15pt}}
$$
So, if $\cA(n)$ is the set of such pairs with $|\la|+s=n$, then the left side of~\eqref{pedod} is $\sum_{n\ge0} \#\cA(n) q^n$.
Similarly, the right side is the generating fuction for 
$$
\cB(n)=\{(\mu,t)\in P_{od}\times \bbN \mid |\mu|+t=n\}.  
$$
So it suffices to construct a bijection $f:\cA(n)\ra\cB(n)$.

To define $f$, suppose $(\la,s)\in\cA(n)$ where $\la=(\la_1,\ldots,\la_l)$.  Let the even parts of $\la$ be $\la_{k+1},\ldots,\la_l$.  
The $f(\la,s)=(\mu,t)$ where
\begin{align*}
    \mu&=(\la_1,\ldots,\la_k,\la_{k+1}-1,\ldots,\la_l-1),\\
    t&=s+l-k.
\end{align*}
Since $\la$ has all parts distinct and had one subtracted from every even part we have $\mu\in P_{od}$.  Furthermore,
$$
|\mu|+t = |\la|-(l-k)+s +(l-k) = |\la|+s =n
$$
So $(\mu,t)\in\cB(n)$ and $f$  is well defined.

To compute $f^{-1}(\mu,t)$ we need to know the number of parts of $\mu$ to which $1$ should be added.  If $f(\la,s)=(\mu,t)$ then this number is controlled by $s$ since if $s>0$ then all parts of $\mu$ should be incremented, whereas if $s=0$ then only the last $t$ should be.
But note that $s=0$ if and only if $t\le \ell(\mu)$.  Indeed if  $s=0$ then, keeping the notation of the previous paragraph and using the fact that $\ell(\la)=\ell(\mu)=l$, we have
$t=0+l-k\le \ell(\mu)$.   And if $s>0$ then $k=0$ so that $t=s+l-0>\ell(\mu)$.
With this in mind, given $(\mu,t)$ with  $l=\ell(\mu)$, we let $m=\min(l,t)$ and define 
$f^{-1}(\mu,t)=(\la,s)$ where
\begin{align*}
    \la&=(\mu_1,\ldots,\mu_{l-m},\mu_{l-m+1}+1,\ldots,\mu_l+1),\\
    s&=t-m.
\end{align*}
It is not hard, using the observation at the beginning of this paragraph, to verify that $f^{-1}$ is well defined and the inverse of $f$, completing the proof.
\eprf

As an example of the map in the previous demonstration, we have
$$
f((9, 5, 4, 2),\ 1) = ((9, 5, 3, 1),\  3).
$$

One can consider the sets $\ol{P}_{xy}^{zw}(n)$ which are the analogues for overpartitions of the 
$P_{xy}^{zw}(n)$. 
Bringmann, Cossaboom and Caig~\cite{BCC:ops} have determined the generating functions of all eight $\pb_{xy}^{zw}$.   The bijective proof of the following result is so like that of Theorem~\ref{Fedod:Thm} that we omit it.
\bth[{\cite[Theorem 1.1 (3)]{BCC:ops}}]
We have

\vs{10pt}

\eqqed{
\Fb_{ed}^{od}(q) = 
\frac{1}{1-q}(-2q;q^2)_\infty - \frac{q}{1-q}(-2q^2;q^2)_\infty.
}
\eth
For comparison with the example after Theorem~\ref{Fedod:Thm}, the map for the result just stated would send
$$
((9, \ol{5}, \ol{4}, 2),\ 1) \mapsto ((9, \ol{5}, \ol{3}, 1),\  3).
$$

\section{Partition triples}
\label{pt}

Guadalupe~\cite{gua:cal} used $q$-series techniques to prove various congruences for the number of partition triples $(\la^{(1)},\la^{(2)},\la^{(3)})$ such that $\la^{(1)}$ and $\la^{(2)}$ have distinct odd parts and $\la^{(3)}$ has all parts divisible by $4$. We will show that one of his results is a special case of a more general theorem which has an easy combinatorial proof.  Furthermore, this result can be generalized further to any prime modulus.

Recall~\eqref{cP} that 
$\cP$ is the set of all partitions
and lete
$$
\cE = \text{ the set of partitions into even parts}.
$$
\bth
\label{cA(n):thm}
Suppose $\cP'\sbe\cP$, $\cE'\sbe\cE$, and 
\beq
\label{cA(n):eq}
\cA(n) = \{(\la^{(1)},\la^{(2)},\la^{(3)}) \mid
\la^{(1)},\la^{(2)}\in\cP',\quad \la^{(3)}\in\cE',\quad |\la^{(1)}|+|\la^{(2)}|+|\la^{(3)}|=n\}.
\eeq
Then
$$
\#\cA(2n+1) \Cong\ 0\ (\Mod 2).
$$
\eth
\bprf
Note that if $(\la^{(1)},\la^{(2)},\la^{(3)})\in\cA(2n+1)$  then $\la^{(1)}\neq\la^{(2)}$.  Indeed, if $\la^{(1)}=\la^{(2)}$ then $|\la^{(1)}|+|\la^{(2)}|$ is even.  
Since $\la^{(3)}\in\cE$, it follows that $|\la^{(1)}|+|\la^{(2)}|+|\la^{(3)}|$ is even. This contradicts the fact that the triple is in $\cA(2n+1)$.

Now define a map $\io:\cA(2n+1)\ra\cA(2n+1)$ by letting
$$
\io(\la^{(1)},\la^{(2)},\la^{(3)})=(\la^{(2)},\la^{(1)},\la^{(3)}).
$$
This is clearly an involution.  But, from the previous paragraph, there are no cycles of length $1$.  Thus $\cA(2n+1)$ decomposes into cycles of length $2$ and so must have even cardinality.
\eprf

As an immediate corollary, we have the following result of Guadalupe.
\bco[{\cite[Theorem 1.1]{gua:cal}}]
Let $\cP'$ be all partitions into distinct odd parts, $\cE'$ be all partitions into parts divisible by $4$, and $\cB(n)$ be the corresponding set of triples as defined by~\eqref{cA(n):eq}.  Then

\vspace{10pt}

\eqqed{\#\cB(2n+1)\Cong 0\ (\Mod 2).}
\eco

The involution used in the proof of Theorem~\ref{cA(n):thm} can be generalized to any prime modulus, $p$, by using the action of the cyclic group of order $p$ on $p$-tuples of partitions by rotation of a tuple.  
For more examples of proofs of congruences using group actions, see the paper of Sagan~\cite{sag:cag}.  This framework is a combinatorial model for 
the following well-known general identity.

\begin{thm}
    For any prime $p$ and any power series $f(q):= \sum_{n=0}^\infty a_n q^n$ with integer coefficients $a_n$, it holds that $f(q)^p:= \sum_{n=0}^\infty b_n q^n$ has the properties that 
    $$
    b_{pn} \equiv a_n \pmod{p},
    $$
    and 
    $$b_{pn+j} \equiv 0 \pmod{p}
    $$ 
    for $j \not\equiv 0 \pmod{p}$.\hqed
\end{thm}

To see why this holds, one considers $f(q)^p$ as generating $p$-tuples of some weighted combinatorial object with counts (possibly negative) given by the coefficients of $f(q)$.  Since $p$ is prime, the only possible periods in the action of cycling the entries of a $p$-tuple are $p$, or 1.  Periods of size 1 correspond to $p$-tuples of identical objects, and other $p$-tuples contribute $0 \pmod{p}$ in total to the coefficients of the power.

We can also give a  characterization of $\#\cB(2n)$ in terms of pairs rather than triples.

\begin{thm}
We have
\beq
\label{cubic}
\#{\mathcal{B}}(2n) \equiv 
\#\{ (\mu^{(1)},\mu^{(2)}) \mid \mu^{(1)}\in\cP,\quad \mu^{(2)}\in\cE,\quad |\mu^{(1)}|+|\mu^{(2)}|=n\}\ \pmod{2}.
\eeq
\end{thm}
\bprf
By a similar argument to that in the demonstration of Theorem~\ref{cA(n):thm}, modulo $2$ it suffices to consider that elements of $\cB(2n)$ where $\la^{(1)}=\la^{(2)}$.
Partitions into distinct odd parts are equinumerous with self-conjugate partitions, so one can replace the first two components with the corresponding self-conjugates $\nu^{(1)}=\nu^{(2)}$.  As for $\la^{(3)}$, since its parts are divisible by $4$, we can replace it by a pair $\nu^{(3)}=\nu^{(4)}\in\cE$ obtained by replacing each part $k$ of $\la^{(3)}$ by a part of $k/2$ in each of the corresponding $\nu$'s.  Thus, we have a bijection up to paritity with pairs $(\nu^{(1)},\nu^{(3)})$ where $\nu^{(1)}$ is self-conjugate, $\nu^{(3)}\in\cE$, and $|\nu^{(1)}|+|\nu^{(3)}|=n$.  Finally, conjugation of the first component is an involution for the pairs on the right side of~\eqref{cubic} whose fixed points are exactly the pairs $(\nu^{(1)}, \nu^{(3)})$ previously constructed.  This finishes the proof.
\eprf

The set on the right side of equation~\eqref{cubic} 
is in bijection with the set of two-colored partitions of $n$ where the odd parts are always red, but the even parts can be either red or blue.  These partitions are the well-known \emph{cubic partitions} of $n$, so named because they arise in the study of Ramanujan's cubic continued fraction.  They were introduced by Chan~\cite{cha:rcc} and named by Kim~\cite{kim:acc}.

\section{Parity-restricted tableaux shapes}
\label{prt}

Partity restrictions also arise in restricting the shapes of standard Young tableaux associated with certain lattice paths.  We begin with the well-known case of Dyck paths.

We will consider lattice paths which are sequences of points $P=p_0,p_1,\ldots,p_n$ in $\bbZ^2$.
We will always take  $p_0=(0,0)$ and $n$ is called the {\em length} of $P$.  The {\em steps} of the path are the vectors $S_i$ from $p_{i-1}$ to $p_i$ for $i\in[n]$ which will be written with square brackets to distinguish them from points in the plane.  Since all our lattice paths  begin at $(0,0)$, we will often write the path as a sequence of steps $P=S_1,S_2,\ldots,S_n$.  

An {\em up step} is $U=[1,1]$; and a {\em down step} is $D=[1,-1]$.  The path $P$ is a {\em Dyck path of semilength $n$} if the following hold
\ben
\item $P$ ends at $(2n,0)$.
\item $P$ consists of $U$ steps and $D$ steps and never goes below the $x$-axis.
\een
Note that a Dyck path of semilength $n$ must have $n$ up steps and $n$ down steps.
Let 
$$
\cD(n) = \{P \mid \text{$P$ is a Dyck path of semilength $n$}\}.
$$
It is well known that 
$$
\#\cD(n)=C_n,
$$
the $n$th Catalan number.

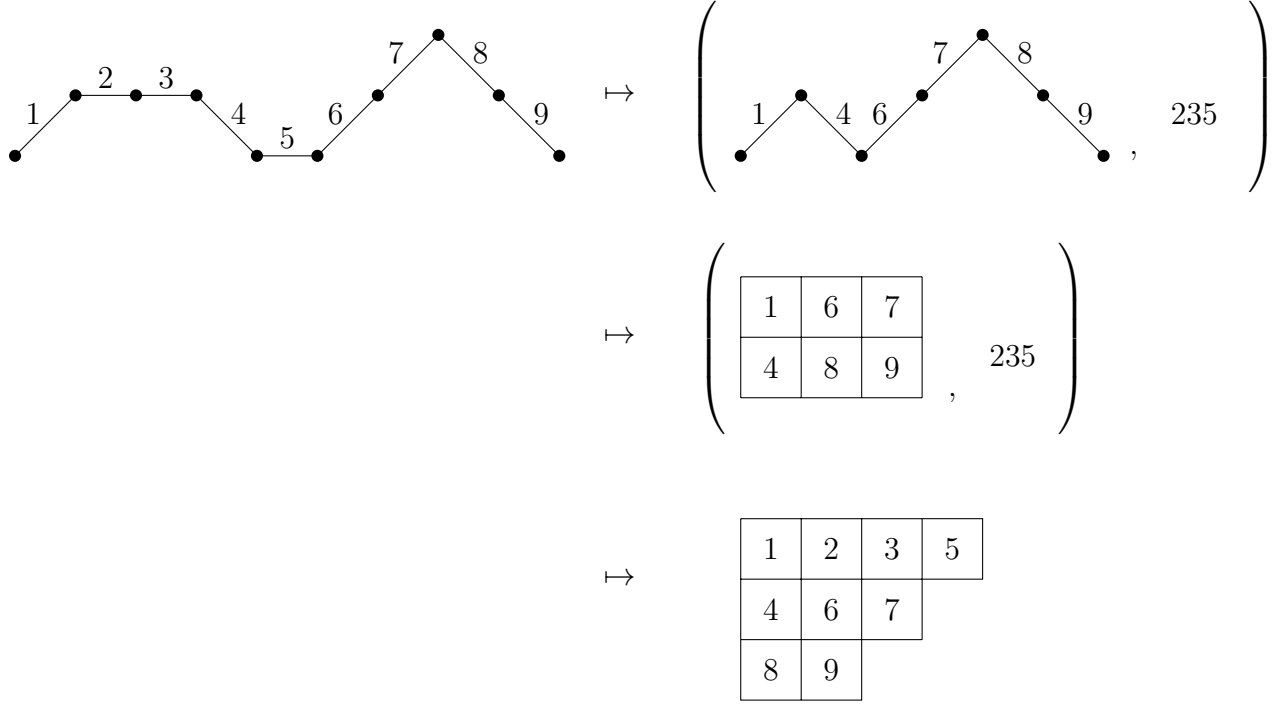
\begin{figure}
    \centering
\begin{tikzpicture}[scale=.8]
\fill(0,0) circle(.1);
\fill(1,1) circle(.1);
\fill(2,1) circle(.1);
\fill(3,1) circle(.1);
\fill(4,0) circle(.1);
\fill(5,0) circle(.1);
\fill(6,1) circle(.1);
\fill(7,2) circle(.1);
\fill(8,1) circle(.1);
\fill(9,0) circle(.1);
\draw (0,0)--(1,1)--(3,1)--(4,0)--(5,0)--(7,2)--(9,0);
\draw(.3,.7) node{$1$};
\draw(1.5,1.3) node{$2$};
\draw(2.5,1.3) node{$3$};
\draw(3.7,.7) node{$4$};
\draw(4.5,.3) node{$5$};
\draw(5.3,.7) node{$6$};
\draw(6.3,1.7) node{$7$};
\draw(7.7,1.7) node{$8$};
\draw(8.7,.7) node{$9$};
\begin{scope}[shift={(11,0)}]
\draw(-1,1) node{$\mapsto$};
\draw(5,1) node{$\left(\rule{0pt}{40pt}\hspace{200pt}\right)$};
\fill(1,0) circle(.1);
\fill(2,1) circle(.1);
\fill(3,0) circle(.1);
\fill(4,1) circle(.1);
\fill(5,2) circle(.1);
\fill(6,1) circle(.1);
\fill(7,0) circle(.1);
\draw (1,0)--(2,1)--(3,0)--(5,2)--(7,0);
\draw(1.3,.7) node{$1$};
\draw(2.7,.7) node{$4$};
\draw(3.3,.7) node{$6$};
\draw(4.3,1.7) node{$7$};
\draw(5.7,1.7) node{$8$};
\draw(6.7,.7) node{$9$};
\draw(7.5,0) node{,};
\draw(8.5,.7) node{$235$};
\end{scope}
\begin{scope}[shift={(11,-4)}]
\draw(-1,1) node{$\mapsto$};
\draw(3.5,1) node{$\left(\rule{0pt}{40pt}\hspace{125pt}\right)$};
\draw (1,0) grid (4,2);
\draw(1.5,1.5) node{$1$};
\draw(2.5,1.5) node{$6$};
\draw(3.5,1.5) node{$7$};
\draw(1.5,.5) node{$4$};
\draw(2.5,.5) node{$8$};
\draw(3.5,.5) node{$9$};
\draw(4.5,0) node{,};
\draw(5.5,.7) node{$235$};
\end{scope}
\begin{scope}[shift={(11,-8)}]
\draw(-1,1) node{$\mapsto$};
\draw (1,0) grid (4,2);
\draw(1.5,1.5) node{$1$};
\draw(2.5,1.5) node{$2$};
\draw(3.5,1.5) node{$3$};
\draw(4.5,1.5) node{$5$};
\draw(1.5,.5) node{$4$};
\draw(2.5,.5) node{$6$};
\draw(3.5,.5) node{$7$};
\draw(1.5,-.5) node{$8$};
\draw(2.5,-.5) node{$9$};
\draw (4,2)--(5,2)--(5,1)--(4,1) (1,0)--(1,-1)--(3,-1)--(3,0) (2,0)--(2,-1)
;
\end{scope}
\end{tikzpicture}
    \caption{An example of the map $m:\cM(9)\ra\cT(9)$}
    \label{m:fig}
\end{figure}

If $\la\ptn n$ is a partition viewed as the set of cells in its  Young diagram then a {\em standard Young tabeau (SYT) of shape $\la$} is a bijection $T:\la\ra[n]$ such that rows and columns increase.
We let
$$
T_{i,j} =\text{ the element in row $i$ and column $j$ of $T$}.
$$
The array in the third line of Fibure~\ref{m:fig} is an SYT, $T$, of shape $(4,3,2)$ with $T_{2,3}=7$.
For more information about SYT, including definitions of any undefined terms used in what follows, see that texts of Sagan~\cite{sag:sg,sag:aoc} or Stanley~\cite{sta:ec2}.
We let
$$
\syt(\la) = \{T \mid \text{$T$ is an SYT of shape $\la$}\},
$$
and
$$
f^\la = \#\syt(\la).
$$

The following result and its proof are well known.  But we include them since they will be useful in the sequel.
\bth
For  $n\ge0$ we have
$$
\#\cD(n) = f^{(n,n)}.
$$
\eth
\bprf
We will define a bijection 
\beq
\label{d:def}
d:\cD(n)\ra\syt(n,n).
\eeq
Given 
$P=S_1 S_2\ldots S_{2n}\in\cD(n)$ we let $d(P)=T$ where
\beq
\label{dcases}
\begin{cases}
 \text{$k$ is in the first row of $T$} & \text{if $S_k=U$},\\ 
 \text{$k$ is in the second row of $T$} & \text{if $S_k=D$},
\end{cases}
\eeq
and the elements in each row are listed in increasing order.
Then $T$ has shape $(n,n)$ since $P$ has $n$ up steps and $n$ down steps, has increasing rows by definition, and increasing columns because $P$ stays weakly about the $x$-axis.  Constructing $d^{-1}$ is left as an exercise for the reader.
\eprf

A {\em horizontal step} is $H=[1,0]$.
A {\em Motzkin path $P$ of length $n$} satisfies the following.
\ben
\item $P$ ends at $(n,0)$.
\item $P$ consists of $U$, $D$ and $H$ steps and never goes below the $x$-axis.
\een
We use the notation
$$
\cM(n) = \{P \mid \text{$P$ is a Motzkin path of length $n$}\}.
$$
The {\em Motzkin numbers} are
$$
M_n = \#\cM(n).
$$
Matsakis and Vendervelde~\cite{MV:mib} gave a bijective proof of the following result.  But it is simpler to use the Robinson-Schendsted map.  For the demonstration, we define the {\em restriction} of $P\in\cM$ to be the Dyck path $P\downarrow$ obtained by removing all the horizontal steps of $P$ and pasting together what remains, keeping the relative order of the steps.  In Figure~\ref{m:fig}, the restriction of the Motzkin path on the left is the Dyck path which is the first component of the first pair on the right.  Such restrictions have played a role in the recent work of Sagan and Sundaram on ordered set partitions~\cite{SS:osp}.

\bth[{\cite[Proposition 4]{MV:mib}}]
For $n\ge0$ we have
$$
\#\cM(n) = \sum_{\scriptsize\stackunder{\la\ptn n}{\ell(\la)\le 3}} f^\la.
$$
\eth
\bprf
Let
$$
\cT(n) = \{T \mid \text{$T\in \syt(\la)$ where $\la\ptn n$ and $\ell(\la)\le3$}\}.
$$
To give a bijection $m:\cM(n)\ra\cT(n)$, take $P\in\cM(n)$ and label its steps with the elements of $[n]$ left to right.  An example of the construction of the map is given in Figure~\ref{m:fig}.  Map $P$ to the pair
$(P\downarrow,L)$ where the labeling of $P$ is preserved in $P\downarrow$, and $L$ is an increasing list of the horizontal steps.
Next, use~\eqref{dcases} to construct a two-row tableau $T$ from $P\downarrow$, again respecting the labeling.  Finally, define $m(P)=U$ where $U$ is the result of inserting $L$ into $T$ using the Robinson-Schensted map.

We claim that $U\in\cT(n)$.  Clearly $U$ is a SYT with $n$ elements by construction.  To check that $U$ has at most $3$ rows note that, since $L$ is increasing, its insertion will cause a horizontal border strip to be added to $T$ from left to right.  Since $T$ has shape $(r,r)$ for some $r$, this forces $U$ to have shape $(r+s,r,t)$ for some $s,t$ with $s+t=\#L$.

To construct the inverse, suppose we are given an SYT, $U$, of shape $\la=(p,r,q)$ where we permit parts equal to $0$.  Use the inverse of Robinson-Schensted to remove all cells of $U$ not in the subshape $(r,r)$, starting with the rightmost cell and working left.  The result will be a pair $(V,L)$ where $V$ has shape $(r,r)$ and $L$ is increasing.  Construct a Dyck path $Q=d^{-1}(V)$ where $d$ is the bijection from~\eqref{d:def} and the labels of $V$ are used on $Q$.  Finally, we obtain a Motzkin path $P$ by inserting horizontal steps in $Q$ labeled by the elements in $L$ so that the labels on $P$ read $1,\ldots,n$ left to right.  This is a step-by-step reversal of the algorithm for the map $m$ and so its inverse.
\eprf

A {\em Riordan path} is a Motzkin path having no horizontal steps on the $x$-axis.  The number of such paths of length $n$ is given by the Riordan number $R_n$.  Let
$$
\cR(n) =\{P \mid \text{$P$ is a Riordan path of length $n$}\}.
$$
We will also be interested in the subsets
$$
\cR(k,m) = \{P  \mid \text{$P$ is a Riordan path with $k$ up steps and $m$ horizontal steps}\}.
$$
Clearly $\cR(n) =\uplus_{2k+m=n}\ \cR(k,m)$.  Hemmer, Straub, and Westrem~\cite{HSW:nic} gave a bijective proof of the following result.
\bth[{\cite[Proposition 1.2]{HSW:nic}}]
For $k,m\ge0$ we have

\vs{10pt}

\eqqed{
\#\cR(k,m) = f^{(k,k,1^m)}.
}
\eth

Hemmer et al.\ also found a description of $\#\cR(n)$ in terms of SYT with parity restricted row lengths, although the proof was not bijective.  To state the next result, call  $\la=(\la_1,\ldots,\la_l)$ {\em nonnegative} if we permit $\la_i =0$ for some $i$.
\bth[{\cite[Theorem 4.5]{HSW:nic}}]
Consider the set of nonnegative partitions
$$
\La(n) =\{ \la=(\la_1,\la_2,\la_3)\ptn n \mid \la_1 \Cong \la_2 \Cong \la_3\ (\Mod 2)\}.
$$
For all $n\ge 0$,

\vs{10pt}

\eqqed{
\#\cR(n) = \sum_{\la\in\La(n)} f^\la.
}
\eth

Combining the previous two results one obtains
$$
\sum_{\la\in\La(n)} f^\la = \sum_{2k+m=n} f^{(k,k,1^m)}.
$$
A refinement of this equation was given in~\cite{HSW:eki}.
Indeed, in the following result that if $n\Cong k\ (\Mod 2)$ then $\La(n,k)\sbe\La(n)$.
\bth[{\cite[Theorem 1.3]{HSW:eki}}]
\label{LnkThm}
Consider the set of nonnegative partitions
$$
\La(n,k) = \{ \la=(\la_1,k,\la_3)\ptn n \mid k \Cong \la_3\ (\Mod 2)\}.
$$
For all $n\ge0$ and $k\ge1$, we have

\vs{10pt}

\eqqed{
\sum_{\la\in\La(n,k)} f^{\la} = f^{(k,k,1^{n-2k})} + f^{(k+1,k+1,1^{n-2k-2})}.
}
\eth

It would be very interesting to find a bijective proof of this theorem.  It is not hard to come up with appropriate maps using conjugation of tableaux when $k=1$ or $2$.

\section{Remarks and open problems}
\label{rop}


We were not able to give bijective proofs for some of the parity-based partition identities which we studied.  We present some of these results here in the hopes that the reader will be able to find suitable maps.

\subsection{More on multiplicity restrictions}

In addition to Theorem~\ref{n=2(mod4)}, Chern was also able to compare $v_e^o(n,0)$ and $v_e^o(n,2)$ when $n\Cong 0\ (\Mod 4)$.  (It is easy to see that when $n$ is odd then both quantities are $0$.)
In particular, he proved
\beq
\label{chern:eq}
\sum_{n\ge0}\ [v_e^o(n,0)-v_e^o(n,2)]\ q^n = \frac{(-q^4;q^4)_\infty}{(q^4;q^8)_\infty}
\eeq
from which it easily follows that
$$
v_e^o(n,0)\ge v_e^o(n,2)
$$
for all $n$ and, in particular, when $4$ divides $n$.  

The Burson-Eichhorn map mentioned earlier also proves this inequality.  But it would be interesting to find a simpler injection.  The problem with using conjugation as we did for $n\Cong 0\ (\Mod 4)$ is that, when $n$ is divisible by $4$, this map carries elements of $V_e^o(n,0)$ to themselves and similarly for $V_e^o(n,2)$.  Still, one could even hope for a nice bijection from $V_e^o(n,0)$ to the disjoint union of $V_e^o(n,2)$ and some set of pairs of partitions naturally counted by the right-hand side of~\eqref{chern:eq}.

\subsection{More on the third order mock theta function}
\label{mto}

Andrews, Dixit, and Yee~\cite{ADY:par}  found a fourth class of partitions counted by the coefficients $p_\nu(n)$ of $\nu(-q)$.  Let $D(n)$ be the set of all partitions  $\la=(\la_1,\la_2,\ldots,\la_l)\vdash n$ satisfying
\ben
\item[(D1)]  $\la_1 > \la_2 >\ldots >\la_l \ge0$, and
\item[(D2)]  for all $i$: if $\la_i$ is odd then $\la_i<2\la_l$.
\een
\bth[{\cite[Theorem 4.1]{ADY:par}}]
For all $n\ge0$ we have

\vs{10pt}

\eqqed{
p_\nu(n) = \# D(n).
}
\eth

It would be interesting to find a combinatorial proof of this theorem, possibly by constructing a bijection between $D(n)$ and one of the sets $A(n)$, $B(n)$, or $C(n)$ from Section~\ref{tom}.  
We will call the partitions in $A(n)$ of {\em type $A$} and similarly for the other sets.
A hint for how to construct such a map is offered by the following result of Andrews and Yee~\cite{AY:imt}.  They prove that an additional parameter can be included by showing the validity of the following two-variable identity. 
\bth[{\cite[Theorem 1, equation (7)]{AY:imt}}]
We have

\vs{10pt}

\eqqed{
\sum_{k\ge0} \frac{x^k q^{k^2+k}}{(q;q^2)_{k+1}} 
= \sum_{j=0}^\infty q^j (-x q^{j+1};q)_j (-x q^{2j+2};q^2)_\infty .
}
\eth  
The coefficient of $x^k$  on the left-hand side of this identity is the generating function for partitions of type $A$ with $k$ distinct even parts, from $2$ through $2k$, so that odds allowed up to $2k+1$.  Meanwhile, the coefficient of $x^k$ on the right-hand side counts partitions of type $D$ with $k+1$ parts.
Hence \emph{some} bijection must exist which matches these refined sets.

Here we return to the map of Chern mentioned at the end of Section~\ref{tom}.  In~\cite{che:cpi}, he considers the three-parameter formula 
$$\sum_{m \geq 0} \frac{x^m q^{m^2+m} }{(yq;q^2)_{m+1}} 
= \sum_{j \geq 0} (yq)^j (-xq/y;q^2)_j .
$$  
He produces a bijective proof of this identity by equating the sizes of the two sets of bipartitions of $n$, that is, pairs of partitions whose sizes sum to $n$:

\begin{itemize}
\item[] \hspace*{-30pt} ${\mathcal{O}}_{m,k}(n)$: bipartitions  $(\lambda, \pi)$ of  $n$ in which $\lambda$ is a partition consisting of $m$ parts of size $m+1$, and $\pi$ is a partition into $k$ odd parts  each of size at most $2m+1$.
\item[] \hspace*{-30pt} ${\mathcal{DO}}_{m,k}(n)$: bipartitions $(\mu, \nu)$ of $n$ in which $\mu$ is a single part of size $m+k$ and $\nu$ is a partition into distinct odd parts with exactly $m$ parts, each of size at most $2(m+k)-1$.
\end{itemize}

Keeping the notation above, Chern's map $\cO_{m,k}(n)\ra\cD\cO_{m,k}(n)$ takes three steps.
\ben
\item[(C1)] For $1\le i\le \ell(\pi)$, if $\pi_i=2s_i+1$ then append $s_i+1$ and $s_i$ to the $i$th row and $i$th column, respectively, of $\la$ to form a partion $\la'$
\item[(C2)] Let $\mu=\la_1'$ .
\item[(C3)] Removing the first row of $\la'$ will form a self-conjugate partition
$\nu'$.  Use the usual bijection between self-conjugate partitions and partitions into odd parts to form $\nu$.
\een

Our map can be produced from Chern's by the following procedure.

\begin{itemize}
\item Replace $\lambda$ with a partition into distinct even parts $2$ through $2m$ and take their disjoint union with the parts of $\pi$ to make a single partition in $A(n)$.
\item In step (C3), 
form a self-conjugate hook with $2|\mu|+1$ cells and paste $\nu'$ inside. Fill the outer hook cells with $1$'s and the rest with $2$'s to produce  the image partition in $B(n)$.
\end{itemize}

The equivalent parameters preserved by our map after this modification are as follows.  In type $A$ partitions, $m$ is half the size of the largest even part and $k$ is the number of odd parts.  In type $B$ partitions,  $m$ is the size of the Durfee square in the self-conjugate partition of $2$'s, and  $m+k$ is the size of the arm of the outermost hook of $1$'s.


With the specialization $y \rightarrow 1$, Chern's identity along with Andrews and Yee's identity implies the following $q$-series identity:

\begin{thm}\label{ACY} $$\sum_{j \geq 0} q^j (-xq;q^2)_j  = \sum_{j \geq 0}
q^j (-xq^{j+1};q)_j (-xq^{2j+2};q^2)_\infty.$$
\end{thm}

This identity certainly does not refine to the level of individual summands, and it would  be of interest to provide a direct bijective proof, perhaps by finding a map between $D(n)$ and $\cD\cO(n)$, where the latter is the union of the $\cD\cO_{m,k}(n)$ over all possible $m$ and $k$.  Ideally, the identity of Andrews and Yee could be further refined with an additional parameter, which might suggest further subdivisions of the equinumerous sets that could guide the creation of such a bijection.

We have made some initial observations and limited bijections concerning edge cases of a map between $D(n)$ and the other sets, including the case for partitions in $D(n)$ where the smallest part is 0 or 1, and partitions in $D(n)$ with up to three parts, i.e., with small powers of $x$ in Chern's identity or Andrews and Yee's identity.  We can provide further details to the interested reader on request.

\subsection{Distinct versus versus repeated parts redux}

Bringmann, Craig and Nazaroglu~\cite{BCN:abp} studied the asymptotic behavior of the cardinalties of the sets $P_{xy}^{zw}(n)$ studied in Section~\ref{dvr}  They proved the following result, asking for combinatorial demonstrations.
\bth[\cite{BCN:abp}]
For sufficiently large $n$ we have
$$
\hspace*{40pt}
p_{ed}^{od}(n)<
p_{od}^{ed}(n)<
p_{od}^{eu}(n)<
p_{eu}^{od}(n)<
p_{ed}^{ou}(n)<
p_{eu}^{ou}(n)<
p_{ou}^{ed}(n)<
p_{ou}^{eu}(n).
\hspace*{40pt}\qed
$$
\eth
Ballantine and Welch~\cite{BW:cpi} were able to prove five of these inequalities, or weaker ones implied by transitivity, using injections.   Then Fan and Xia~\cite{FX:cpi} gave an injective map for the fourth inequality above.  It would be pleasing to provide such map for the remaining case of the fourth inequality.

\nocite{*}
\bibliographystyle{alpha}

\end{document}